\documentclass[Bernoulli,preprint]{imsart}

\usepackage{amsthm,amsmath}
\RequirePackage[colorlinks,citecolor=blue,urlcolor=blue]{hyperref}

\usepackage[round]{natbib}
\usepackage{graphicx}
\startlocaldefs
\endlocaldefs

\begin{document}

\begin{frontmatter}

\title{The Circular Mat\'ern Covariance Function and its Link to Markov Random Fields on the Circle}
\runtitle{Circular Mat\'ern Covariance Functions}

\begin{aug}

\author{\fnms{Chunfeng} \snm{Huang}\corref{}\ead[label=e1]{huang48@indiana.edu}},
\author{\fnms{Ao} \snm{Li} \ead[label=e2]{liao@alumni.iu.edu}}

\address{Department of Statistics, Indiana University\\ Bloomington, IN 47408, USA }
 \affiliation{Indiana University}

\printead{e1,e2}
\phantom{E-mail:\ }



\end{aug}

\begin{abstract}
The link between Gaussian random fields and Markov random fields is well established based on a stochastic partial differential equation in Euclidean spaces, where the Mat\'ern covariance functions are essential.  However, the Mat\'ern covariance functions are not always positive definite on circles and spheres. In this manuscript, we focus on the extension of this link to circles, and show that the link between Gaussian random fields and Markov random fields on circles is valid based on the circular Mat\'ern covariance function instead. First, we show that this circular Mat\'ern  function is the covariance of the stationary solution to the stochastic differential equation on the circle with a formally defined white noise space measure. Then, for the corresponding conditional autoregressive model, we derive a closed form formula for its covariance function.  Together with a closed form formula for the circular Mat\'ern covariance function, the link between these two random fields can be established explicitly.  Additionally, it is known that the estimator of the mean is not consistent on circles, we provide an equivalent Gaussian measure explanation for this non-ergodicity issue.

\end{abstract}



\end{frontmatter}

\section{Introduction.} Gaussian random fields and Markov random fields are two important sub-areas in spatial statistics \citep{Cressie1993, CressieandWikle2011}. They each possess unique and different models and methods. Gaussian random fields are studied extensively in geostatistics, where covariance functions play the essential role in modeling spatial dependency. On the other hand, Markov random fields focus on conditional distributions and precision matrices.  The link between these two random fields has been established in the celebrated work by \citet{Lindgrenetal2011}. For a Gaussian random field when its covariance function is Mat\'ern, its link to the Markov random field is based on a stochastic partial differential equation (SPDE). \citet{Lindgrenetal2011} establish this link mainly in Euclidean spaces, and briefly discuss the extension to other manifolds, especially circles and spheres. However, some recent works have shown that these Mat\'ern covariance functions are not always positive definite on circles and spheres \citep{Huangetal2011, Gneiting2013}. This poses the question on the extension from Euclidean spaces to other manifolds. In this manuscript, we focus on circles and show that the extension to circles is valid when the circular Mat\'errn covariance functions are used instead. \\

We first set to find the solution of the SPDE in \citet{Lindgrenetal2011} on circles, where a proper white noise space measure needs to be defined.  Note that the white noise is not an ordinary random process, but a generalized random process \citep{Ito1953, GelfandandVilenkin1964}. One common way is to view it as a generalized derivative of a Brownian motion \citep{Kuo1996}. However, the Brownian motion on circles constructed in  \cite{Levy1959} appears to be problematic, and is shown to be a regular Euclidean Brownian motion on the half circle, but the exact mirror image on the other half  \citep{HuangandLi2021}. With this observation, \citet{HuangandLi2021} formally introduce a white noise measure on the circle through the dual Sobolev spaces. Based on this development, we derive the solution to the SPDE, and obtain the corresponding covariance function in Section 2. This matches a type of covariance functions introduced in \cite{GuinnessandFuentes2016}, which they name it a circular Mat\'ern covariance function.  \\

Given the development of such Gaussian random field with the circular Mat\'ern covariance function on the circle, we proceed to find its connection to Markov random fields in Section 3. First, for a circular Mat\'ern covariance function of order 1, a corresponding conditional autoregressive model (CAR) is constructed. For this CAR model, we invert the precision matrix to obtain its covariance function. We then show how to find a closed form formula for this covariance function. Note that a closed form formula for the circular Mat\'ern covariance function of order $1$ is readily available \citep{GuinnessandFuentes2016}, and these two formulae are shown to exactly match each other. Therefore, the link between these two random fields is exact and explicit. We continue with the circular Mat\'ern covariance function of order $2$ and build a corresponding CAR model. For this CAR model, a closed form formula is also derived and is shown to approximate the circular Mat\'ern covariance function of order $2$. It is clear from our findings that the extension of the link between Gaussian random fields and Markov random fields to circles is valid when the circular Mat\'ern covariance functions are used.  \\

In this manuscript, it is worth noting that we make a few additional discoveries. The white noise measure sheds light on a non-ergodicity issue on circles. \cite{Lauritzen1973} discovered that ergodicity and Gaussian cannot coexist on circles and spheres \citep{Schaffrin1993}. In Remark 2.2 , we provide an explanation based on the white noise space measure \citep{HuangandLi2021} and the equivalent Gaussian measures. In the computational front, we find an alternative and simpler way to obtain the closed form expression for the circular Mat\'ern covariance functions compared to \cite{GuinnessandFuentes2016} in Appendix B. \\

Circular spaces may be the most rudimentary manifold. We hope the results and methods developed in this note pave the way for further understanding of the Gaussian random fields, Markov random fields, and their connection in other manifolds. \\

\section{Circular Mat\'ern covariance function and SPDE.} Mat\'ern covariance functions are popularly used in modeling Gaussian random fields in spatial statistics \citep{Cressie1993, Stein1999, CressieandWikle2011}. In this manuscript, we use the terms random processes and random fields exchangeably. In Euclidean spaces,  a process  $\{ X(t), t \in R^d\}$ is assumed to have the Mat\'ern covariance function \citep{Stein1999}, if its covariance adapts the form
\[
\mbox{cov} (X(t), X(s)) = \frac{\sigma^2}{2^{\nu-1} \Gamma(\nu)} (\kappa \|t-s\|)^{\nu} K_{\nu} (\kappa \| t - s \|), \quad s, t \in R^d, 
\]
where $K_{\nu}(\cdot)$ is the modified Bessel function of the second kind of order $\nu > 0$, $\kappa > 0$ and $\sigma^2$ are parameters, and $\| \cdot \|$ is the Euclidean distance. This function is shown to be a covariance function of a stationary solution of the following stochastic partial differential equation (SPDE) \citep{Whittle1954, Whittle1963, Besag1981, Lindgrenetal2011, Vergaraetal2022}  
\begin{equation} \label{eq:eq1}
(\kappa^2 - \Delta)^{\alpha/2} X(t) = W(t), \quad \alpha = \nu + d/2, \quad \nu > 0, \quad \kappa > 0, \quad t \in R^d,
\end{equation}
where $\Delta$ is the Laplacian operator in the Euclidean space $R^d$, $W(t)$ is the white noise process, $(\kappa^2 - \Delta)^{\alpha/2}$ is a pseudo-differential operator. Based on this, \cite{Lindgrenetal2011}  showed the connection between Gaussian random fields and Markov random fields, and extended this link to other manifolds, for example, circles and spheres. However, noted in  \cite{Huangetal2011} and \cite{Gneiting2013}, the Mat\'ern covariance functions are not positive definite on circles and spheres when $\nu > 1/2$. Therefore, such extension to circles proposed in \cite{Lindgrenetal2011} is called in question, and becomes the focus of this manuscript. To investigate this extension, we first study the SPDE ({\ref{eq:eq1}) on the circles when $t \in S$, where $S$ is a unit circle, that is, the following SPDE
\begin{equation} \label{eq:eqSPDEcircle}
(\kappa^2-\Delta)^{\alpha/2} X(t) = W(t), \quad t \in S, \quad \alpha > \frac{1}{2}.
\end{equation} 
The resulting covariance function will be valid on circles, and provide the basis for building the link between the Gaussian random fields and Markov random fields. \\

First, a white noise process on circles needs some care. Based on the concept of the generalized random processes \citep{GelfandandVilenkin1964}, \cite{HuangandLi2021} have formally introduced a white noise measure on the circle, and will be used in this Section. Then, we follow the S-transform in \cite{Si2012} to obtain the solution to the SPDE (\ref{eq:eqSPDEcircle}) on the circle. We start with the regular $L^2$ space on the circle $S$, where
\[
L^2(S) = \{ f(t), t \in S, \int_S |f(t)|^2 dt < \infty \},
\]
with the inner product $(f,g) = \int_S f(t) \bar{g}(t) dt$, and the associated norm $\| f \|^2 = (f, f)$.  Then, a white noise measure is the triple \citep{HuangandLi2021}
\[
(H_{-1}(S), \mathcal{B}, \mu),
\]
where $H_{-1}(S)$ is the Sobolev space of index $(-1)$ on the circle, $\mathcal{B}$ is the Borel $\sigma$-algebra on $H_{-1}(S)$, $\mu(H_{-1}(S))=1$ is the white noise measure.  The derivation relies on the Gel'fand triple
\[
H_{1,0}(S) \subset L_0^2(S) \subset H_{-1}(S),
\]
where $H_{1,0}(S)$ is the dual space of $H_{-1}(S)$, which is a Sobolev space of index $1$ with one extra condition that 
\begin{equation} \label{eq:eqextracondition}
\int_S f(t) d t = 0, \quad f \in H_{1,0}(S),
\end{equation}
and $L_0^2(S)$ is the $L^2(S)$ space with this same extra condition (\ref{eq:eqextracondition}). Given this white noise measure space, we follow \cite[Section 4.3]{Si2012} and use S-transform to find the solution of the SPDE (\ref{eq:eqSPDEcircle}). Consider the S-transform:
\[
(S \phi) (\xi) = C(\xi) \int_{H_{-1}(S)} e^{(\omega, \xi)} \phi(\omega) d \mu(\omega), \quad \xi \in H_{1,0}(S), 
\]
where 
\[
C(\xi) = e^{-\frac{1}{2} \| \xi \|^2}
\]
is the characteristic functional of the white noise measure and $\phi(\cdot) \in L^2(H_{-1}(S), \mathcal{B}, \mu)$. Applying S-transform to both sides of the SPDE (\ref{eq:eqSPDEcircle}), we obtain
\begin{equation} \label{eq:ode}
(\kappa^2 - \Delta)^{\alpha/2} U(t, \xi) = \xi(t),
\end{equation}
where
\[
U(t, \xi) = (S X(t)) (\xi), \quad \xi \in H_{1,0}(S).
\]
This becomes an ordinary differential equation with the pseudo-differential operator defined through Fourier transform \citep{Samkoetal1992} 
\[
\{ \mathcal{F} (\kappa^2 - \Delta)^{\alpha/2} U \} (k) = (\kappa^2 + (2 \pi k)^2 )^{\alpha/2} \{ \mathcal{F} U \} (k), \quad k=0, \pm 1, \ldots.
\]
For $\xi \in H_{1,0}(S)$, its Fourier expansion is
\[
\xi = \sum_{k=-\infty}^{\infty} \xi_k e^{i 2 \pi k t},
\]
that is, 
\[
\{ \mathcal{F} \xi \} (k) = \xi_k,
\]
Therefore,
\[
\{ \mathcal{F} U \}(k) = \frac{\xi_k}{(\kappa^2 + (2 \pi k)^2)^{\alpha/2}}, \quad k=0, \pm 1, \ldots.
\]
Solving the ordinary differential equation (\ref{eq:ode}) leads to
\[
U(t, \xi) = \sum_{k = -\infty}^{\infty} \frac{\xi_k}{(\kappa^2 + (2 \pi k)^2)^{\alpha/2}} e^{i 2 \pi k t}.
\]
This can be written in $L_2^0(S)$ inner product as
\[
U(t, x) = (G(t,u), \xi(u)),
\]
where
\[
G(t,u) = \sum_k \frac{1}{(\kappa^2 + (2 \pi k)^2)^{\alpha/2}} e^{i 2 \pi k (t-u)}.
\]
This is the Green's function of the ordinary differential equation (\ref{eq:ode}). By taking the inverse S-transform, we obtain the solution to the SPDE (\ref{eq:eqSPDEcircle}):
\[
X(t) = (G(t,u), W(u)).
\]
For this random process, the covariance function \citep[see][Lemma 6]{HuangandLi2021} is
\begin{equation} \label{eq:eq2}
\mbox{cov} (X(t), X(s)) = (G(t,u), G(s,u)) = \sum_k \frac{e^{i 2 \pi k(t-s)}}{(\kappa^2 + (2 \pi k)^2)^{\alpha}}.
\end{equation}
This is the covariance function of the stationary solution of the SPDE (\ref{eq:eqSPDEcircle}). For the range of $\alpha$, since $G(t, u) \in L_0^2(S)$, it is clear that $\alpha > 1/2$.\\

{\bf Remark 2.1.} The covariance (\ref{eq:eq2}) matches a type of covariance function introduced in \citet[equation (7)]{GuinnessandFuentes2016}, which they  name it the circular Mat\'ern covariance function. Our derivation here shows that such covariance functions can be directly obtained through SPDE (\ref{eq:eqSPDEcircle}).  In this manuscript, we term the random field with the covariance function (\ref{eq:eq2}) as the circular Mat\'ern random field. \\

{\bf Remark 2.2. on Non-ergodicity.} \citet[Section 8.3 and Section 10.2]{Lauritzen1973}  shows that the estimators of the mean and the covariance function are not consistent for a homogenous (i.e., stationary) process on the sphere. It was also stated in \cite{Schaffrin1993} that the homogeneous processes on spheres which are both Gaussian and ergodic do not exist. Certainly, there is the same problem for random processes on circles. In particular, if a stationary process $X(t)$ on a unit circle is assumed to have the mean $\mbox{E} X(t) = \mu$ and the covariance function
\[
\mbox{cov} (X(t), X(s)) = C(t-s) = a_0 + \sum_{n=1}^{\infty} a_n \cos n (t-s).
\]
This expansion can be found in \cite{Schoenberg1942} or \cite{Huangetal2016}. Even with the ability of observing the entire process on the circle, the estimators
\[
\hat{\mu} = \frac{1}{2 \pi} \int_S X(t) d t,
\]
and
\[
\hat{C} (h)  = \frac{1}{ (2 \pi)^2} \int \int_{S \times S, d(t,s) = h} X(t) X(s) d t d s,
\]
where $d(t,s)$ is angular distance between $t$ and $s$, are not consistent because their variances do not vanish. The white noise space measure developed in \cite{HuangandLi2021} can help explain this. Given $X(t)$ and an arbitrary uncorrelated random variable
\[
X_0 \sim N(\mu_0, \sigma_0^2),
\]
Now, consider another random process
\[
Y(t) = X(t) + X_0.
\]
The characteristic functionals for $X(t)$ and $Y(t)$ rely on $(X(t), \xi)$ and $(Y(t), \xi)$, where $\xi \in H_{1,0}(S)$. Note that $\xi$ satisfies the extra condition (\ref{eq:eqextracondition}), which results in $(X_0, \xi) = 0$. Therefore, 
\[
(Y(t), \xi) = (X(t) + X_0, \xi) = (X(t), \xi).
\]
That is, the characteristic functionals of $X(t)$ and $Y(t)$ are the same, and the Gaussian measures on these two processes are equivalent. This implies that there cannot be consistent estimators of either $\mu_0$ or $\sigma_0^2$. Similar phenomena can be also found in \citet[Chapter 3]{Wahba1990a}. Note that this extra condition (\ref{eq:eqextracondition}) leads naturally to the Brownian bridge, instead of Brownian motion on the circle \citep{HuangandLi2021}. \\

\section{CAR models and the Link between two random fields.}  In this section, we study the conditional autoregressive (CAR) models on the circles, and establish the link between circular Mat\'ern  random fields and Markov random fields. For a CAR model, one can obtain its precision matrix and the corresponding covariance function by inverting the precision matrix. We discover a closed form formula for this covariance function. Together with the closed form expression for the circular Mat\'ern covariance function, the link can be shown explicitly. \\

{\bf CAR model for $\alpha=1$.} \\

We start with the circular Mat\'ern random field when $\alpha=1$ with the following covariance function ($\alpha=1$ in equation (\ref{eq:eq2})), 
\begin{equation} \label{eq:eq4}
\mbox{cov} (X(t), X(s)) = \sum_{k=-\infty}^{\infty} \frac{e^{i 2 \pi k(t-s)}}{\kappa^2 + (2 \pi k)^2}, \quad t, s \in [0,1].
\end{equation}
To build the link, consider a CAR model on equally-spaced grids on a unit circle:
\[
\{ Z(\theta_k) \}, \quad \theta_k = \frac{2 \pi}{n} k, \quad k=1, 2, \ldots, n,
\]
where we assume the conditional distribution is Gaussian
\begin{equation} \label{eq:car2}
Z(\theta_k) | Z(\theta_{-k}) \sim N(a Z(\theta_{k-1}) + aZ(\theta_{k+1}), \sigma^2), \quad a > 0, \quad \theta_{n+k} = \theta_k.
\end{equation}
Then, the joint distribution of $Z = (Z(\theta_1), \ldots, Z(\theta_n))^T$ can be shown to be  \citep{Besag1974, Cressie1993} 
\[
Z \sim N(0, \sigma^2 (I-M_1)^{-1}),
\]
where
\[
M_1 = \left( \begin{array}{cccccc} 0 & a & 0 & \cdots & \cdots & a \\ 
a & 0 & a & \cdots & \cdots & 0 \\
0 & a & 0 & \cdots & \cdots & 0 \\
\vdots & \vdots & \vdots & \ddots & \vdots & \vdots \\
0 & 0 & \cdots & \cdots & 0 & a \\
a & 0 & \cdots & \cdots & a & 0 \end{array} \right).
\]
This is a circulant matrix, and has a spectral decomposition
\[
M_1 = P \Lambda P^*,
\]
where
\[
\Lambda = \mbox{diag} \{ a(e^{i \frac{2 \pi}{n} (k-1)} + e^{ - i \frac{2\pi}{n} (k-1)}) \}_{k=1, \ldots, n},
\]
and 
\[
P = \frac{1}{\sqrt{n}} \{ e^{- i \frac{2\pi}{n} (k_1-1)(k_2-1)} \}_{n \times n}, \quad k_1, k_2=1,\ldots,n,
\]
and $P^*$ is its Hermitian.  While the matrix $M_1$ and the precision matrix $\frac{1}{\sigma^2} (I-M_1)$ are sparse, the covariance matrix $\sigma^2(I-M_1)^{-1}$ is not sparse. By the spectral decomposition the covariance matrix $\sigma^2(I-M_1)^{-1} =\sigma^2 P(I-\Lambda)^{-1} P^*$ and we have
\begin{equation} \label{eq:eq3}
\mbox{cov} (Z(\theta_{k_1}), Z(\theta_{k_2})) = \frac{\sigma^2}{n} \sum_{k=1}^n \frac{e^{i \frac{2\pi}{n} (k_2-k_1) (k-1)}}{1-a(e^{i \frac{2\pi}{n} (k-1)} + e^{-i \frac{2\pi}{n} (k-1)})}, \quad k_1, k_2 = 1, \ldots, n.
\end{equation}

The similarity and difference between these two covariance functions (\ref{eq:eq4}) and (\ref{eq:eq3}) are quite pronounced in two ways: (I) the summation in (\ref{eq:eq3}) is through $n$, and the summation in (\ref{eq:eq4}) is through $\infty$; (II) if one conducts a Taylor expansion of the cosine function in the denominator in (\ref{eq:eq3}) with
\begin{equation} \label{eq:eqstar}
1- a(e^{i \frac{2\pi}{n}(k_1)} +e^{-i \frac{2\pi}{n} (k-1)}) = 1- 2 a \cos \frac{2\pi}{n} (k-1) \approx (1-2a) + \frac{a}{n^2} ( 2 \pi (k-1))^2,
\end{equation}
this mimics the denominator in equation (\ref{eq:eq4}). In Euclidean spaces, this is how \cite{Besag1981} proposed the approximation of the Mat\'ern covariance function for the CAR model in $R^2$, where the integration extends from $[0,\pi]$ to $[0, \infty)$, and the Taylor expansion of the cosine function is also used. \cite{Besag1981}'s approach serves as the basis in \cite{Lindgrenetal2011}  to establish the fundamental connection between Gaussian random fields and Markov random fields in Euclidean spaces. It is still possible to extend such approximation to circles, see Remark 3.2 below. However, we discover that there are closed form expressions for both summations (\ref{eq:eq4}) and (\ref{eq:eq3}). This makes the connection explicit and the approximation used in \cite{Besag1981} is not necessary. \\

For circular Mat\'ern covariance  (\ref{eq:eq4}), one can use Equation (1.445.2) in \citet{GradshteynandRyzhik1994} \citep{GuinnessandFuentes2016} and obtain the closed form formula in hyperbolic functions
\begin{equation} \label{eq:eq5}
\mbox{cov} (X(t), X(s)) = \frac{1}{2 \kappa \sinh \frac{\kappa}{2}} \cosh \left\{ \kappa (|t-s| - \frac{1}{2}) \right\}.
\end{equation}
For CAR equation (\ref{eq:eq3}), we can factor the denominator and derive a closed form expression (details can be found in Appendix A)
\begin{equation} \label{eq:eq6}
\mbox{cov} (Z(\theta_{k_1}), Z(\theta_{k_2})) =   \frac{\sigma^2}{ \tanh (\log \beta) \sinh \frac{n \log \beta}{2}} \cosh \left\{ n \log \beta \cdot (\frac{ |k_1-k_2|}{n} - \frac{1}{2} ) \right\},
\end{equation}
where
\[
\beta = \frac{1+\sqrt{1-4a^2}}{2a}.
\]

{\bf Remark 3.1.} While one can see the potential connection between two covariance functions (\ref{eq:eq4}) and (\ref{eq:eq3}), the closed form expressions (\ref{eq:eq5}) and (\ref{eq:eq6}) make their link much more transparent. In particular, given a circular Mat\'ern random field with the covariance function (\ref{eq:eq4}) and the parameter $\kappa$, we can build a CAR model (\ref{eq:car2}) with an arbitrary $n$,
\[
a= \frac{1}{2 \cosh \frac{\kappa}{n}} \quad \mbox{and} \quad \sigma^2 =  \frac{\tanh \frac{\kappa}{n}}{2 \kappa}.
\]
Then, $n \log \beta = \kappa$, and for this CAR model, the covariance 
\begin{equation} \label{eq:eq7}
\mbox{cov} (Z(\theta_{k_1}), Z(\theta_{k_2})) = \frac{1}{2 \kappa \sinh \frac{\kappa}{2}} \cosh \left\{ \kappa ( \frac{|k_1-k_2|}{n} - \frac{1}{2}) \right\}.
\end{equation}
This is exactly the same as the circular Mat\'ern covariance function (\ref{eq:eq5}), where $\frac{|k_1-k_2|}{n}$ is in the place of $|t-s|$.  That is, given the circular Mat\'ern random field with covariance (\ref{eq:eq3}), one can build a CAR model which yields the same covariance structure. Reversely, given a CAR model (\ref{eq:car2}) with $a, \sigma^2$ and $n$, one can find a corresponding circular Mat\'ern random field with $\alpha=1, \kappa=n \cosh^{-1} \frac{1}{2a}$ and variance $2n \sigma^2 \log \beta/\sqrt{1-4a^2}$. Therefore, the equivalence between the Gaussian random fields  and the Markov random fields when $\alpha=1$ is established. \\

{\bf Remark 3.2.}  While the two closed form formulae provide the explicit exact link between the two random fields when $\alpha=1$, one can also follow Besag (1981) to build the CAR model to approximate the circular Mat\'ern random field. For example, by Taylor Expansion of the denominator (\ref{eq:eqstar}), we can approximate equation (\ref{eq:eq3}),  
\begin{eqnarray*}
\mbox{cov}(Z(\theta_{k_1}, Z(\theta_{k_2}))  \approx  \sigma^2 \frac{n}{a} \sum_{k=1}^n \frac{e^{i \frac{2\pi}{n} (k_2-k_1)(k-1)}}{\frac{n^2(1-2a)}{a} + 4 \pi^2 (k-1)^2}.
\end{eqnarray*}
If we follow the similar approach in \citet{Besag1981}, and match
\[
\kappa^2 = \frac{n^2 (1-2a)}{a}.
\]
We obtain
\[
a = \frac{n^2}{\kappa^2 + 2n^2}.
\]
Compare this with the previous match by Taylor expansion and assuming $n$ is large,
\[
a = \frac{1}{2 \cosh \frac{\kappa}{n}} \approx \frac{1}{2 (1 + \frac{1}{2} \frac{\kappa^2}{n^2})} = \frac{1}{2 + \frac{\kappa^2}{n^2}} = \frac{n^2}{\kappa^2 + 2n^2}.
\]
Similarly, $\beta \approx e^{\kappa/n}$.  Note that, \cite{Besag1981} approach was introduced in Euclidean spaces, its application on the circle  will provide an approximation. It is clear that the exact match in Remark 3.1 will be preferred, and this approach is not necessary. \\

{\bf CAR model for $\alpha=2$.} \\

Now, let us consider the circular Mat\'ern random field when $\alpha=2$ with the following covariance function,
\begin{equation} \label{eq:eq7}
\mbox{cov} (X(t), X(s)) = \sum_{k=-\infty}^{\infty} \frac{e^{i 2 \pi k(t-s)}}{(\kappa^2 + (2 \pi k)^2)^2}, \quad t, s \in [0,1].
\end{equation}
We build a CAR model through the convolution of the CAR model (\ref{eq:car2}):
\begin{equation} \label{eq:eqCAR4model}
Z(\theta_k) | Z(\theta_{-k}) \sim N(a_1 Z(\theta_{k-1}) + a_1 Z (\theta_{k+1}) + a_2  Z(\theta_{k-2}) + a_2 Z(\theta_{k+2}), \sigma^2), 
\end{equation}
where
\[
a_1 = \frac{2a}{2a^2+1}, \quad a_2 =- \frac{a^2}{2a^2+1}.
\]
This convolution approach follows \citet{Lindgrenetal2011}. For this CAR model,  the covariance matrix is $\sigma^2 (I-M_2)^{-1}$, where
\[
M_2 = \left[ \begin{array}{cccccccc} 0 & a_1 & a_2 & 0 & \cdots & 0 & a_2 & a_1 \\
a_1 & 0 & a_1  & a_2 & \cdots & 0 & 0 & a_2 \\
a_2 & a_1 & 0 & a_1 & \cdots & 0 & 0 & 0 \\
\vdots & \vdots & \vdots & \vdots & \vdots & \vdots & \vdots &\vdots \\
a_1 & a_2 & 0 & 0 & \cdots & a_2 & a_1 & 0 \end{array} \right].
\]
Use the property of the circulant matrix, we can derive the corresponding covariance function
\begin{equation} \label{eq:eq9}
\mbox{cov} (Z(\theta_{k_1}), Z(\theta_{k_2})) = \sigma^2 \frac{2a^2+1}{n} \sum_{k=1}^n \frac{e^{i \frac{2\pi}{n} (k_2-k_1)(k-1)}}{ \left( 1- a (e^{i \frac{2\pi}{n} (k-1)} + e^{- i \frac{2\pi}{n} (k-1)}) \right)^2}.
\end{equation}

Similar to the case of $\alpha=1$, one can clearly see the difference and similarity between the two covariance functions (\ref{eq:eq7}) and (\ref{eq:eq9}). Again, we discover that there are closed form expressions for both summations, which will help us establish the link. First,  equation (\ref{eq:eq7}) yields (see \citealp{GuinnessandFuentes2016}, Appendix B, and remark 3.3 below)
\begin{eqnarray} \label{eq:eq8}
&\,& \mbox{cov}(X(t), X(s)) =  \frac{\sinh \frac{\kappa}{2} + \frac{\kappa}{2} \cosh \frac{\kappa}{2}}{4 \kappa^3 \sinh^2 \frac{\kappa}{2}}  \cosh \left\{ \kappa (|t-s|-\frac{1}{2}) \right\} \nonumber \\
&\,& \quad \quad \quad -  \frac{|t-s|-\frac{1}{2}}{4 \kappa^2 \sinh \frac{\kappa}{2}} \sinh \left\{ \kappa (|t-s|-\frac{1}{2}) \right\}.
\end{eqnarray}

{\bf Remark 3.3.}  \cite{GuinnessandFuentes2016} show the summation (\ref{eq:eq8}) through the differential relationship with respect to $|t-s|$ in (\ref{eq:eq4}). We find an alternative and simpler way by taking the derivative with respect to the parameter $\kappa$ instead. The details are provided in Appendix B. \\

As for the CAR model's covariance (\ref{eq:eq9}), we derive its closed form expression (see Appendix A): 
\begin{eqnarray} \label{eq:eqCAR4}
&\,&  \mbox{cov} (Z(\theta_{k_1}), Z(\theta_{k_2})) =  n \sigma^2 (2a^2+1)  \nonumber \\
&\,& \hspace{-0.6cm} \times  \left[ \frac{\frac{1}{2} \cosh \frac{n \log \beta}{2} + \frac{1}{n} \coth (\log \beta) \sinh \frac{n \log \beta}{2} }{\tanh^2 (\log \beta) \sinh^2 \frac{n \log \beta}{2}}  \cosh \left\{ (n \log \beta) ( \frac{|k_1-k_2|}{n} - \frac{1}{2} )\right\} \right.  \nonumber \\
&\,& \hspace{-0.95cm} \left. -  \frac{1}{\tanh^2 (\log \beta) \sinh \frac{n \log \beta}{2}}( \frac{|k_1-k_2|}{n} - \frac{1}{2}) \sinh \left\{ (n \log \beta) (\frac{|k_1-k_2|}{n} - \frac{1}{2}) \right\} \right],
\end{eqnarray}
where
\[
\beta = \frac{1+\sqrt{1-4a^2}}{2a}.
\]
Then, let
\[
\sigma^2 = \frac{\sinh^2 ( \log \beta)}{2 n^3 \log^2 \beta (1+2 \cosh^2 (\log \beta))},
\]
the covariance (\ref{eq:eqCAR4}) becomes
\begin{eqnarray} \label{eq:eq11}
&\,& \mbox{cov} (Z(\theta_{k_1}), Z(\theta_{k_2}))  \nonumber \\
&\,& =  \frac{ \frac{n \log \beta}{2} \cosh \frac{n \log \beta}{2} + \log \beta \coth (\log \beta) \sinh \frac{n \log \beta}{2}}{4 (n \log \beta)^3 \sinh^2 \frac{n \log \beta}{2}}  \cosh \left\{ (n \log \beta ) ( \frac{|k_1-k_2|}{n} - \frac{1}{2})  \right\} \nonumber \\
&\,& \quad -  \frac{1}{4 (n \log \beta)^2 \sinh \frac{n \log \beta}{2}} ( \frac{|k_1-k_2|}{n} - \frac{1}{2}) \sinh \left\{ (n \log \beta) ( \frac{|k_1-k_2|}{n} - \frac{1}{2}) \right\}.
\end{eqnarray}

Therefore, given a circular Mat\'ern covariance function (\ref{eq:eq8}), we can build the CAR model (\ref{eq:eqCAR4model})  with $a = \frac{1}{2 \cosh \frac{\kappa}{n}}$, or $n \log \beta  = \kappa$, which results in the following covariance function:
\begin{eqnarray} \label{eq:eqCAR4kappa}
&\,& \mbox{cov} (Z(\theta_{k_1}), Z(\theta_{k_2}))  \nonumber \\
&\,& =  \frac{1}{4 \kappa^3 \sinh^2 \frac{\kappa}{2}} \left( \frac{\kappa}{2} \cosh \frac{\kappa}{2} + \frac{\kappa}{n} \coth \frac{\kappa}{n} \sinh \frac{\kappa}{2} \right) \cosh \left\{ \kappa ( \frac{|k_1-k_2|}{n} - \frac{1}{2}) \right\}  \nonumber \\
&\,& \quad -  \frac{1}{4 \kappa^2 \sinh \frac{\kappa}{2}} ( \frac{|k_1-k_2|}{n} - \frac{1}{2}) \sinh \left\{ \kappa ( \frac{|k_1-k_2|}{n} - \frac{1}{2}) \right\}.
\end{eqnarray}
This time, these two functions, (\ref{eq:eq8}) and (\ref{eq:eqCAR4kappa}) are not exactly the same. The CAR model covariance function (\ref{eq:eqCAR4kappa}) differs from the circular Mat\'ern covariance function (\ref{eq:eq8}) with the term
\[
\frac{\kappa}{n} \coth \frac{\kappa}{n}
\]
in front of $\sinh \frac{\kappa}{2}$, instead of just $1$. That is, this CAR model (\ref{eq:eqCAR4model})  only approximates the circular Mat\'ern covariance function when $\alpha=2$. For this approximation, when $n$ increases, 
\[
\frac{\kappa}{n} \coth \frac{\kappa}{n} \approx \frac{\kappa}{n} \left( \frac{n}{\kappa} + \frac{\kappa}{3n} \right) = 1 + \frac{\kappa^2}{3 n^2}.
\]
That is, for a fixed $\kappa$, one can increase $n$ to make this as close to $1$ as possible. Therefore, given $\kappa$, one can construct a CAR model (\ref{eq:eqCAR4model}) to approximate (\ref{eq:eq8})  with an $n$ such that $\kappa \ll n$, and
\[
a = \frac{1}{2 \cosh \frac{\kappa}{n}}, \quad \mbox{and} \quad \sigma^2 = \frac{\sinh^2 \frac{\kappa}{n}}{2 n \kappa^2 (1+2 \cosh^2 \frac{\kappa}{n})}.
\]

In Figure 1, we show two plots of the circular Mat\'ern correlation function when $\alpha=2$ with the same $\kappa=10$, and the left panel is with $n=10$, and right panel is with $n=50$. The solid line is the correlation function of the circular Mat\'ern correlation function from (\ref{eq:eq8}), and the dotted line is the correlation function from (\ref{eq:eqCAR4kappa}) of the corresponding CAR model (\ref{eq:eqCAR4model}). It is clear that when $n$ is relatively small with respect to $\kappa$, one might see slight difference in correlation functions (left panel). However, this difference diminishes rapidly with a larger $n$ (right panel).   \\ 

\begin{figure}[h]
\begin{center}
\includegraphics[width=10cm,height=6cm]{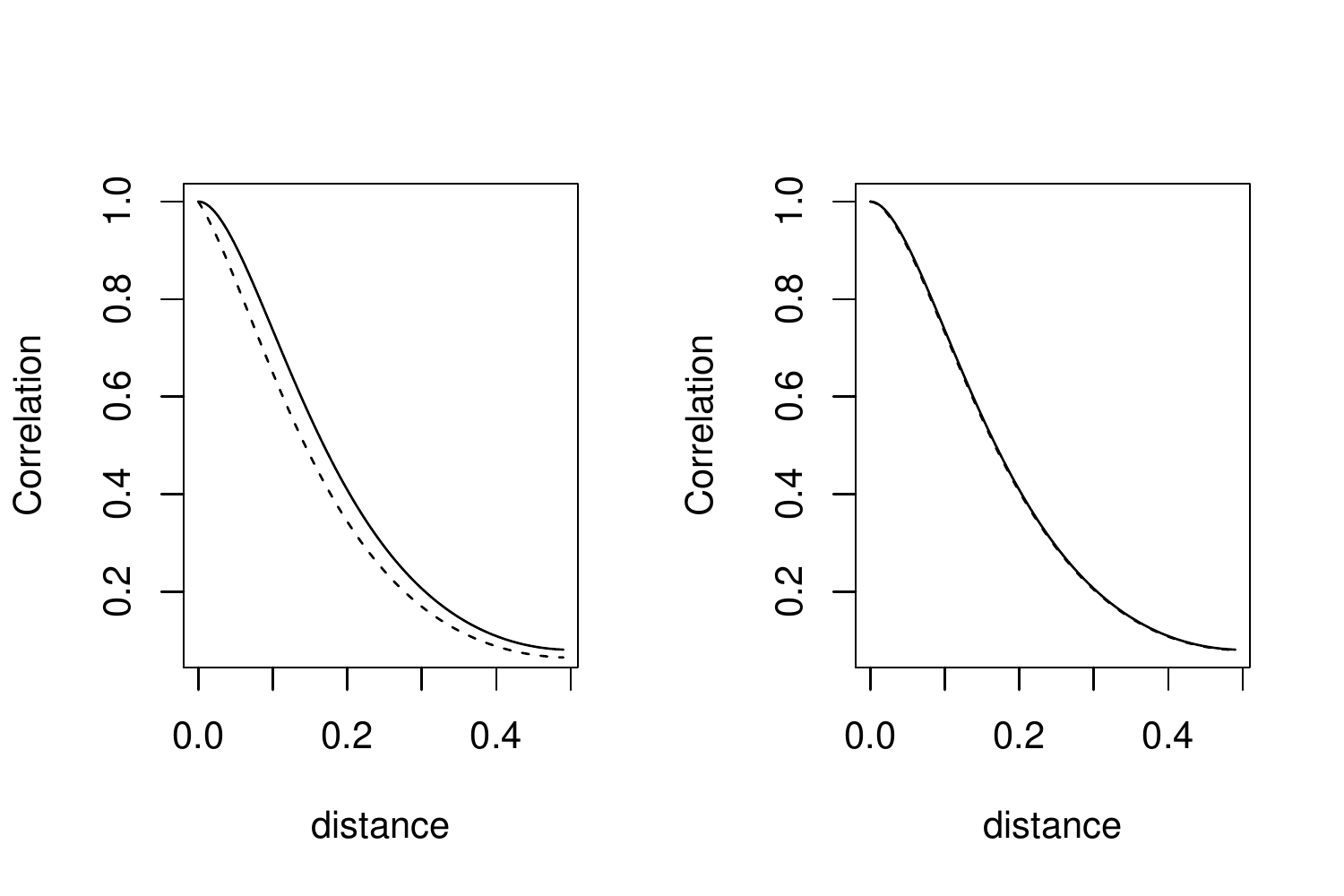}
\caption{The solid line is the circular Mat\'ern correlation function with $\alpha=2$, and the dotted line is the correlation function for the corresponding CAR model. Left panel is with $\kappa=10, n=10$, and right panel is with $\kappa=10, n=50$.}
\end{center}
\end{figure}


{\bf Remark 3.4.} Similar to {Remark 3.1}, for a circular Mat\'ern random field with $\alpha=2$, we build a CAR model (\ref{eq:eqCAR4model}). The closed form expressions (\ref{eq:eq8}) and (\ref{eq:eqCAR4kappa}) for both random fields reveal the striking similarity. This time, the match is not exact, but an approximation. This is different from {Remark 3.1}. Nevertheless,  given $\kappa$, one can build a CAR model to approximate the circular Mat\'ern random field with $\alpha=2$. The link between these two random fields is clear. Combining this remark with {Remark. 3.1}, we see that the extension of \cite{Lindgrenetal2011} to circles are valid when the circular Mat\'ern covariance function is used instead. \\

{\bf Remark 3.5.}  In Appendix A, we show how to derive  (\ref{eq:eqCAR4}) , where two approaches are presented. One is similar to the CAR model (\ref{eq:car2})  derivation, which is lengthy. The other way is to directly take the derivative with respect to $a$ in (\ref{eq:eq6}), and is much simpler. More details are in Appendix A. \\

{\bf Remark 3.6.} For the circular Mat\'ern covariance functions with higher orders one can follow our approach in Appendix B or \cite{GuinnessandFuentes2016} for derivation. For example, equation (\ref{eq:eq14}) in Appendix B shows how to obtain the circular Mat\'ern covariance function when $\alpha=3$. Similarly, for the corresponding CAR models, the general formulae can be obtained, for example, see equation (\ref{eq:eq13}) in Appendix A.  These formulae will help establish the explicit link,  but appear to be lengthy. We hope that what we learn in this note paves the way for building the link between Gaussian random fields and Markov random fields on spheres and other manifolds.   \\

{\bf Appendix A}.\\

In this Appendix, we show how to obtain closed form formulae for CAR models covariance functions, in particular, equations (\ref{eq:eq6}) from (\ref{eq:eq3}) and (\ref{eq:eqCAR4}) from (\ref{eq:eq9}). First, for $m=1,2, \ldots,$ let
\[
\phi_m(\theta) = \frac{1}{n} \sum_{k=1}^n \frac{e^{i \frac{2\pi}{n} (k-1) \theta}}{ \left( 1- a (e^{i \frac{2\pi}{n} (k-1)} + e^{-i \frac{2\pi}{n} (k-1)}) \right)^m}, \quad \theta = \frac{|k_1-k_2|}{n} \in [0,1].
\]
These are CAR covariance functions in this manuscript with a scaling factor. We simplify the notation from $\frac{|k_1-k_2|}{n}$ to $\theta$. For the CAR model (\ref{eq:car2}), the covariance function (\ref{eq:eq3}) will be
\[
\mbox{cov} (Z(\theta_{k_1}), Z(\theta_{k_2})) = \sigma^2 \phi_1(\theta).
\]
For  this $\phi_1(\theta)$, we note that the denominator can be factored into the product 
\begin{equation} \label{eq:eq12}
1-a(e^{i \frac{2\pi}{n} (k-1)} + e^{-i \frac{2\pi}{n} (k-1)}) = (x-be^{i \frac{2\pi}{n} (k-1)})(x-be^{- i \frac{2\pi}{n} (k-1)}),
\end{equation}
where
\[
x = \frac{\sqrt{1+2a}- \sqrt{1-2a}}{2}, \quad \mbox{and,} \quad b=\frac{\sqrt{1+2a}+\sqrt{1-2a}}{2}.
\]
Let 
\[
\beta = \frac{b}{x} = \frac{1+\sqrt{1-4a^2}}{2a},
\]
then we have the expansion
\begin{eqnarray*}
&\,& \frac{1}{x-b e^{i \frac{2\pi}{n} (k-1)}} = \frac{1}{x} \frac{1}{1-\beta e^{i \frac{2\pi}{n}(k-1)}} = \frac{1}{x(1-\beta^n)} \sum_{t=1}^{n} \beta^{t-1} e^{i \frac{2\pi}{n} (k-1)(t-1)},
\end{eqnarray*}
and
\[
\frac{1}{x-be^{-i\frac{2\pi}{n} (k-1)}} = \frac{1}{x(1-\beta^n)} \sum_{t=1}^n \beta^{t-1} e^{- \frac{2\pi}{n} (k-1)(t-1)}.
\]
Together, we have 
\begin{eqnarray*}
&\,& \phi_1(\theta) = \frac{1}{nx^2(1-\beta^n)^2} \sum_{k=1}^n \sum_{t_1=1}^n \sum_{t_2=1}^n e^{i \frac{2\pi}{n} (k_2-k_1) (k-1)} \beta^{t_1+t_2-2} e^{i \frac{2\pi}{n} (k-1) (t_1-t_2)} \\
&\,& = \frac{1}{n x^2 (1-\beta^n)^2} \sum_{t_1=1}^n \sum_{t_2=1}^n \beta^{t_1+t_2-2} \sum_{k=1}^n e^{i \frac{2\pi}{n} (k-1)[ (t_1-t_2)-(k_1-k_2)]}.
\end{eqnarray*}
The inside summation has the property that 
\[
\sum_{k=1}^n e^{i \frac{2\pi}{n} (k-1)[(t_1-t_2)-(k_1-k_2)]} = \left\{ \begin{array}{l} n, \quad \mbox{if } (t_1-t_2) - (k_1-k_2) \quad \mbox{mod} \quad n = 0, \\ 0, \quad \mbox{otherwise}. \end{array} \right.
\]
It is clear that $\phi_1(\cdot)$ is an even function in $(k_1-k_2)$, here, we assume $k_1-k_2 \ge 0$.  Then the summation is $n$ when $t_1$ is from $1$ through $k_1-k_2$, and the corresponding $t_2 = n-(k_1-k_2)+t_1$, and when $t_1$ is from $(k_1-k_2+1)$ through $n$, the corresponding $t_2 = t_1-(k_1-k_2)$. In the other words, the summation is not zero when
\[
(t_1, t_2) \in \{ (1, n-(k_1-k_2)+1), \ldots, (k_1-k_2,n) \} \cup \{ (k_1-k_2+1, 1), \ldots, (n, n-(k_1-k_2))\}.
\]
Therefore, the summation is composed of two parts:
\begin{eqnarray*}
\phi_1(\theta) = \frac{1}{x^2(1-\beta^n)^2} \left( \sum_{t_1=1}^{k_1-k_2} \beta^{t_1  + (n-(k_1-k_2)+t_1)-2} + \sum_{t_1=k_1-k_2+1}^n \beta^{t_1+(t_1-(k_1-k_2))-2} \right),
\end{eqnarray*}
which leads to
\[
\phi_1(\theta) = \frac{\beta^{n-(k_1-k_2)} + \beta^{k_1-k_2}}{x^2 (1-\beta^n) (1-\beta^2)}.
\]
With the following equalities:
\begin{eqnarray*}
&\,& \beta = \frac{1+\sqrt{1-4a^2}}{2a}, \quad a = \frac{1}{2 \cosh (\log \beta)},\\
&\,& x^2(1-\beta^2) = - \sqrt{1-4a^2}= - \tanh (\log \beta) , \\
&\,& \frac{\beta^{n-(k_1-k_2)} + \beta^{k_1-k_2}}{1-\beta^n} = - \frac{\beta^{\frac{n}{2}}}{\beta^n-1} \left( \beta^{\frac{n}{2} - (k_1-k_2)} + \beta^{(k_1-k_2) - \frac{n}{2}} \right) \\
&\,& \quad = - \frac{ \cosh \left\{ n \log \beta \cdot (\frac{k_1-k_2}{n} - \frac{1}{2}) \right\}}{\sinh \frac{n \log \beta }{2}}.
\end{eqnarray*}
We obtain
\[
\phi_1(\theta) = \frac{1}{\tanh (\log \beta) \sinh \frac{n \log \beta}{2}} \cosh \left\{ n \log \beta (\theta - \frac{1}{2})\right\}.
\]
and arrive at equation (\ref{eq:eq6}). \\

Now for $\phi_2(\theta)$, we show how to obtain equation (\ref{eq:eqCAR4}). First, there is a differential relationship between $\phi_1(\theta)$ and $\phi_2(\theta)$ with respect to $a$:
\[
\phi_2(\theta) = \phi_1(\theta) + a \frac{d}{da} \phi_1(\theta).
\]
In fact,  by rescaling the denominator of $\phi_1(\theta)$ and taking the derivative with respect to $a$, we have
\begin{eqnarray*}
&\,& \phi_1(\theta) = \frac{1}{an} \sum_{k=1}^n \frac{e^{i 2 \pi (k-1) \theta}}{\frac{1}{a} - (e^{i \frac{2\pi}{n} (k-1)} + e^{-i \frac{2\pi}{n} (k-1)})}, \\
&\,& \frac{d}{d a} \phi_1(\theta) = - \frac{1}{a^2 n} \sum_{k=1}^n \frac{e^{i 2 \pi (k-1) \theta}}{\frac{1}{a} - (e^{i \frac{2\pi}{n} (k-1)} + e^{-i \frac{2\pi}{n} (k-1)})} \\
&\,& \quad \quad + \frac{1}{an} \sum_{k=1}^n \frac{1}{a^2} \frac{e^{i 2 \pi (k-1) \theta}}{ \left( \frac{1}{a} - (e^{i \frac{2\pi}{n}(k-1)} + e^{-i \frac{2\pi}{n} (k-1)}) \right)^2 } \\
&\,& = - \frac{1}{a} \phi_1(\theta) + \frac{1}{a} \phi_2(\theta).
\end{eqnarray*}
Based on this, we can compute and obtain  
\begin{eqnarray*}
&\,& \phi_2(\theta) =n \frac{\frac{1}{2} \cosh \frac{n \log \beta}{2} + \frac{1}{n} \coth (\log \beta) \sinh \frac{ n \log \beta}{2}}{\tanh^2( \log \beta) \sinh^2 \frac{n \log \beta}{2}} \cosh \left\{ n \log \beta (\theta - \frac{1}{2} )  \right\} \\
&\,& \quad - \frac{n}{\tanh^2 (\log \beta) \sinh \frac{n \log \beta}{2}} (\theta - \frac{1}{2}) \sinh \left\{ n \log \beta (\theta - \frac{1}{2}) \right\}. \\
\end{eqnarray*}

In our manuscript, the equation (\ref{eq:eq9}) of the CAR model (\ref{eq:eqCAR4model}) is
\[
\mbox{cov}(Z(\theta_{k_1}), Z(\theta_{k_2})) = \sigma^2 (2a^2+1) \phi_2(\theta).
\]
Then, we obtain equation (\ref{eq:eqCAR4}). \\

{\bf Remark A.1.} Follow this direction of taking the derivative with respect to $a$, one can obtain the closed form formulae for $\phi_m(\theta)$ when $m$ increases. For example, when $m=3$, we have
\begin{equation} \label{eq:eq13}
\phi_3(\theta) = \phi_2(\theta) + \frac{a}{2} \frac{d}{da} \phi_2(\theta).
\end{equation}

\vspace{0.4cm}

{\bf Remark A.2.}  Here, we provide an alternative way to derive the closed form for $\phi_2(\theta)$. One can follow the approach in $\phi_1(\theta)$ and factor the denominator. First,
\begin{eqnarray*}
&\,& \phi_2(\theta) = \frac{1}{n} \sum_{k=1}^n \frac{e^{i \frac{2\pi}{n} (k_2-k_1) (k-1)}}{\left( 1- a(e^{i \frac{2\pi}{n} (k-1)} + e^{-i \frac{2\pi}{n} (k-1)}) \right)^2} \\
&\,& = \frac{1}{n} \sum_{k=1}^n \frac{e^{i \frac{2\pi}{n} (k_2-k_1)(k-1)}}{(x-be^{i \frac{2\pi}{n}(k-1)})^2 (x-be^{-i \frac{2\pi}{n} (k-1)})^2},
\end{eqnarray*}
where $x$ and $b$ are the same in equation (\ref{eq:eq12}), and the same $\beta = \frac{b}{x}$. Then, 
\begin{eqnarray*}
&\,& \phi_2(\theta) = \frac{1}{n x^4(1-\beta^n)^4} \sum_{k=1}^n e^{i \frac{2\pi}{n} (k_2-k_1)(k-1)} \\
&\,& \quad \times \sum_{t_1=1}^n \sum_{t_2=1}^n \sum_{t_3=1}^n \sum_{t_4=1}^n \beta^{t_1+t_2+t_3+t_4-4} e^{i \frac{2\pi}{n} (k-1) ((t_1-t_2) + (t_3-t_4))}
\end{eqnarray*}
Following the same arguments, the last summation is $n$ only when 
\[
(t_1-t_2) + (t_3-t_4) - (k_1-k_2) \quad \mbox{mod} \quad n = 0. 
\]
Then, the summation will be decomposed into eight terms:
\begin{eqnarray*}
&\,& \phi_2(\theta) =  \frac{1}{x^4(1-\beta^n)^4} \\
&\,&  \hspace{-.4in} \times \left[ \sum_{u_1=0, u_1+u_2=k_1-k_2}^{k_1-k_2} \left( \sum_{t_1=1,t_1-t_2=u_1-n}^{u_1 > 0} + \sum_{t_1=u_1+1,t_1-t_2=u_1}^n \right) \left( \sum_{t_3=1,t_3-t_4=u_2-n}^{u_2 > 0} + \sum_{t_3=u_2+1, t_3-t_4=u_2}^n \right) \right. \\ 
&\,& \hspace{-.8in}  + \sum_{u_1=(k_1-k_2)+1, u_1+u_2=k_1-k_2+n}^{n-1} \left( \sum_{t_1=1,t_1-t_2=u_1-n}^{u_1 > 0} + \sum_{t_1=u_1+1,t_1-t_2=u_1}^n \right) \left( \sum_{t_3=1,t_3-t_4=u_2-n}^{u_2 > 0} + \sum_{t_3=u_2+1, t_3-t_4=u_2}^n \right) \\
&\,& \quad   \beta^{t_1+t_2+t_3+t_4-4} \Bigg].\\
\end{eqnarray*}
Carefully going through these eight terms, we can obtain  
\begin{eqnarray*}
&\,& \phi_2(\theta) = \frac{1}{x^4(1-\beta^n)^2 (1-\beta^2)^2} \\
&\,& \times \left( (|k_1-k_2|+1) (\beta^{2n-|k_1-k_2|} + \beta^{|k_1-k_2|}) + (n-1-|k_1-k_2|) (\beta^{n+|k_1-k_2|} + \beta^{n-|k_1-k_2|}) \right. \\
&\,&  \quad \left. + \frac{2 (1-\beta^n)}{1-\beta^2} (\beta^{n-|k_1-k_2|} + \beta^{|k_1-k_2|+2}) \right).
\end{eqnarray*}
Noting that 
\[
\beta^{|k_1-k_2|-\frac{n}{2}} = \cosh \left\{ n \log \beta (\theta - \frac{1}{2}) \right\} + \sinh \left\{ n \log \beta (\theta - \frac{1}{2}) \right\},
\]
and
\[
\beta^{\frac{n}{2} - |k_1-k_2|} = \cosh \left\{ n \log d (\theta-\frac{1}{2})  \right\} - \sinh \left\{ n \log \beta (\theta-\frac{1}{2}) \right\},
\]
and with the algebraic properties of hyperbolic functions, we arrive exactly the same $\phi_2(\theta)$.  This approach follows the derivation in $\phi_1(\theta)$, but it appears to be very lengthy. \\

{\bf Appendix B}. \\

In this Appendix, we show an alternative way to obtain the closed form expression for the circular Mat\'ern covariance function (\ref{eq:eq2}).  Let
\[
\psi_m (\theta) = \sum_{k=-\infty}^{\infty} \frac{e^{i 2 \pi k \theta}}{(\kappa^2 + (2\pi k)^2)^m}, \quad \theta \in [0,1], \quad m=1, 2, \ldots
\]
When $m=1, 2$, these become equations (\ref{eq:eq5}) and (\ref{eq:eq8}), respectively. For $m=1$, we can obtain the summation directly from (1.445.2) in \cite{GradshteynandRyzhik1994}
\[
\psi_1(\theta) = \frac{1}{2 \kappa \sinh \frac{\kappa}{2}} \cosh \left\{ \kappa (\theta- \frac{1}{2}) \right\}.
\]
For $m=2$, one can follow \cite{GuinnessandFuentes2016} that there is a differential relationship with respect to $\theta$:
\[
\kappa^2 \psi_2(\theta) - \frac{d^2}{d \theta^2} \psi_2(\theta) = \psi_1(\theta)
\]
and compute the value of $\psi_2(\frac{1}{2})$ to obtain the closed form formula.\\

Here, we provide an alternative way. Consider taking the derivative of $\psi_1(\theta)$ with respect to $\kappa$: 
\begin{eqnarray*}
&\,& \frac{d}{d \kappa} \psi_1(\theta) = - 2 \kappa \sum_{k=-\infty}^{\infty} \frac{e^{i 2 \pi k \theta}}{ (\kappa^2 + (2 \pi k)^2)^2} = - 2 \kappa \psi_2(\theta).
\end{eqnarray*}
That is, we can obtain $\psi_2(\theta)$ more directly:
\begin{eqnarray*}
&\,& \psi_2(\theta) = - \frac{1}{2 \kappa} \cdot \frac{d}{ d \kappa} \psi_1(\theta) = - \frac{1}{2 \kappa} \cdot \frac{d}{d \kappa} \frac{1}{2 \kappa \sinh \frac{\kappa}{2}} \cosh \left\{ \kappa (\theta - \frac{1}{2}) \right\} \\
&\,& = \frac{\sinh \frac{\kappa}{2} + \frac{\kappa}{2} \cosh \frac{\kappa}{2}}{4 \kappa^3 \sinh^2 \frac{\kappa}{2}} \cosh \left\{ \kappa (\theta - \frac{1}{2}) \right\} - \frac{1}{4 \kappa^2 \sinh \frac{\kappa}{2}} (\theta - \frac{1}{2}) \sinh \left\{ \kappa (\theta - \frac{1}{2}) \right\}.
\end{eqnarray*}
We arrive at equation (\ref{eq:eq8}). \\

{\bf Remark B.1.} This approach is an alternative to \cite{GuinnessandFuentes2016}, and appears to be simpler. \\

{\bf Remark B.2.} Similar to {Remark A.1}, we can extend this approach to  obtain, $\psi_m(\theta)$ for general integer $m$. For example, when $m=3$,
\begin{equation} \label{eq:eq14}
\psi_3(\theta)= - \frac{1}{4 \kappa} \frac{d}{d \kappa} \psi_2(\theta).
\end{equation}

\bibliography{mybibfile}

\end{document}